\documentclass{amsart}
\usepackage{amssymb}
\usepackage{graphicx}

\def\wr{\,\operatorname{wr}}
\def\PL{\operatorname{PL}}

\def\supp{\operatorname{supp}}
\def\Z{\mathbf Z}
\def\R{\mathbf R}
\def\inv{{}^{-1}}
\def\str#1{\langle#1\rangle}
\def\av#1{\overline{#1}}
\def\avst#1{\overline{\mathstrut #1}}
\def\cM{\mathcal M}
\def\cN{\mathcal N}
\def\a{\alpha}

\newtheorem{Th}{Theorem}

\newtheorem{Lem}[Th]{Lemma}
\newtheorem{Prop}[Th]{Proposition}
\newtheorem{Claim}[Th]{Claim}

\renewcommand{\le}{\leqslant}
\renewcommand{\ge}{\geqslant}

\begin{document}

\title[Interpreting the arithmetic in Thompson's group $F$]{Interpreting the arithmetic in \\ Thompson's group $F$}
\author[Valery Bardakov and Vladimir Tolstykh]{Valery Bardakov${}^\dagger$ \and Vladimir Tolstykh}
\address{Valery Bardakov\\ Institute of Mathematics\\
Siberian Branch Russian Academy of Science\\
630090 Novosibirsk\\
Russia}
\email{bardakov@math.nsc.ru}
\address{Vladimir Tolstykh\\ Department of Mathematics\\ Yeditepe University\\
34755 Kay\i\c sda\u g\i \\
Istanbul\\
Turkey}
\email{vtolstykh@yeditepe.edu.tr}
\thanks{${}^\dagger$The first author gratefully acknowledges support from T\"UB\.ITAK, the Scientific
and Research Council of Turkey}
\subjclass[2000]{20F65 (03C35, 03C60)}
\maketitle

\begin{abstract}
We prove that the elementary theory of Thompson's group $F$ is hereditarily
undecidable.
\end{abstract}

M.~Sapir asked whether the elementary theory of Thompson's group $F$
is decidable \cite[a part of Question 4.16]{TGr40yrs}.
More specifically, he asked later whether $F$ first-order interprets
the restricted wreath product $\Z \wr \Z.$ It is known that the
elementary theory of the latter group first-order interprets
the arithmetic, and is therefore hereditarily undecidable \cite{Nosk}. Recall that a first-order theory $T$ is said to be
{\it hereditarily undecidable} if every subtheory of
$T$ is undecidable.

The aim of the present note is to show that $F$ indeed first-order
interprets $\Z \wr \Z,$ and hence the elementary theory
of $F$ is hereditarily undecidable.

Let $J$ be any closed interval of $\R$ with dyadic endpoints.
The group $\PL_2(J)$ is the group of all continuous order-preserving
piecewise linear maps from $J$ into itself with finitely
many breakpoints all occurring at dyadic rationals and with the slopes of their linear `components'
equal to powers of two. One of the standard realizations of Thompson's
group $F$ is given by the group $\PL_2([0,1])$
(for the background information on Thompson's group $F$ the reader
may consult \cite{Belk,CFP}.) It is well-known that
$\PL_2([0,1])$ is generated by the following maps $x_0$ and $x_1$:
$$
\alpha x_0=
\begin{cases}
2\alpha,                    &0 \le \alpha \le \frac 14, \\
\a+\frac 14,               &\frac 14 \le \alpha \le \frac 12, \\
\frac 12 \a +\frac 12,    &\frac 12 \le \alpha \le 1,
\end{cases}
\qquad
\alpha x_1=
\begin{cases}
\alpha,                    &0 \le \alpha \le \frac 12, \\
2\alpha-\frac12,                   &\frac 12 \le \alpha \le \frac 58, \\
\a+\frac 18,               &\frac 58 \le \alpha \le \frac 34, \\
\frac 12 \a +\frac 34,     &\frac 34 \le \alpha \le 1.
\end{cases}
$$
Note that we write maps of $\PL_2([0,1])$ on the right.

By the definition,
$$
x_n = x_0^{-n+1} x_1 x_0^{n-1}
$$
where $n \ge 1$ is a natural number.

For any map $f \in \PL_2([0,1])$ we define
\begin{align*}
& \supp(f)=(\text{the support of $f$})=\{\alpha \in [0,1] : \alpha f \ne
\alpha\}.
\end{align*}
Clearly,
\begin{align*}
&\supp(f\inv g f)=\supp(g) f
\end{align*}
for all $f,g \in \PL_2([0,1]).$ It follows that if maps $f,g \in \PL_2([0,1])$
have disjoint supports, then they are commuting.

The restricted wreath product $\Z \wr \Z$ has the following presentation:
$$
\Z \wr \Z \cong \str{a,b\ |\ \{ [b,a^{-n} b a^n] : n \in \Z\}}.
$$

The group $F$ is `riddled' with copies of $\Z \wr \Z,$ as,
for instance, the results from \cite{GS} show.
We shall demonstrate that one such a copy is given by the group $G$ generated by
$a=x_0^2$ and $b=x_1 x_2^{-1}$ (cf. \cite{Clea} where it is
proved that the subgroup generated by $x_0$ and $x_1 x_2 x_1^{-2}$
is also isomorphic to $\Z \wr \Z.$)

\begin{Lem} \label{is_isom}
{\em (i)} The supports $S_k$ of the maps $a^{-k} b a^k=x_0^{-2k} x_{1} x_{2}^{-1} x_0^{2k}$
where $k \in \Z$ are pairwise disjoint open intervals such that
$$
\bigcup_{k \in \Z} \avst S_k =[0,1],
$$
where $\av J$ is the closure of an interval $J$ of $\R;$

{\em (ii)} The group $G$ is isomorphic to $\Z \wr \Z.$
\end{Lem}

\begin{proof}
(i) A straightforward calculation shows that
$$
\supp(x_1x_2^{-1}) =\left(\frac 12, \frac 78\right)=\left(\frac 12, \frac 12 x_0^2\right)
$$
(to calculate the products of elements of $\PL_2([0,1])$ we suggest
the reader use the method of the tree diagrams, see, for instance,
\cite[Section 1.2]{Belk}.)

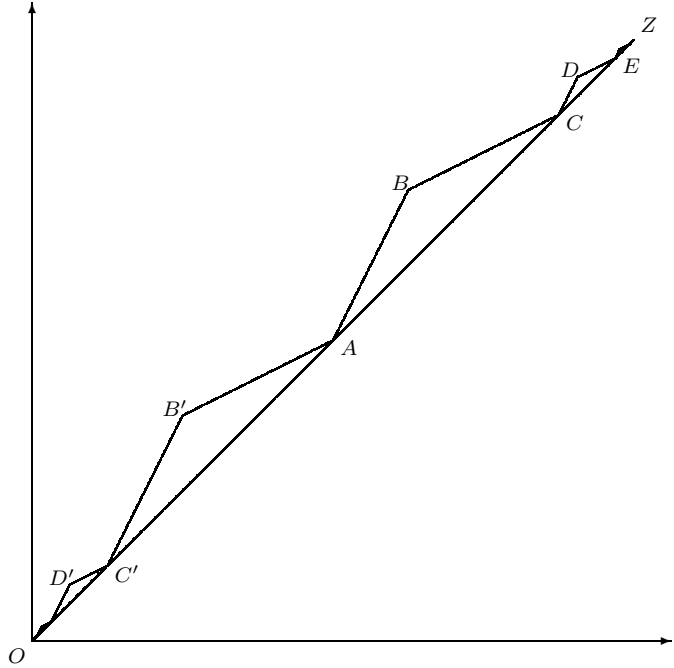
\begin{figure}[h]

\unitlength=1.00mm
\special{em:linewidth 0.4pt}
\linethickness{0.4pt}
\begin{picture}(90,90)
\put(5,5){\vector(1,0){85}}
\put(5,5){\vector(0,1){85}}


\put(3,3){\makebox(0,0)[cc]{\footnotesize $O$}}
\qbezier(5,5)(25,25)(45,45)
\put(46,44){\makebox(0,0)[lc]{\footnotesize $A$}}
\qbezier(45,45)(50,55)(55,65)
\put(54,66){\makebox(0,0)[cc]{\footnotesize $B$}}
\qbezier(55,65)(65,70)(75,75)
\put(76,74){\makebox(0,0)[lc]{\footnotesize $C$}}
\qbezier(75,75)(80,80)(85,85)
\put(87,87){\makebox(0,0)[cc]{\footnotesize $Z$}}


\qbezier(5,5)(40,40)(75,75)
\qbezier(75,75)(76.25,77.5)(77.5,80)
\put(76.5,81){\makebox(0,0)[cc]{\footnotesize $D$}}
\qbezier(77.5,80)(80,81.25)(82.5,82.5)
\put(83.5,81.5){\makebox(0,0)[lc]{\footnotesize $E$}}
\qbezier(82.5,82.5)(83.75,83.75)(85,85)


\qbezier(5,5)(43.75,43.75)(82.5,82.5)
\qbezier(82.5,82.5)(82.8125,83.125)(83.125,83.75)
\qbezier(83.125,83.75)(83.75,84.0625)(84.375,84.375)
\qbezier(84.375,84.375)(84.6875,84.6875)(85,85)


\qbezier(5,5)(10,10)(15,15)
\put(16,14){\makebox(0,0)[lc]{\footnotesize $C'$}}
\qbezier(15,15)(20,25)(25,35)
\put(24,36){\makebox(0,0)[cc]{\footnotesize $B'$}}
\qbezier(25,35)(35,40)(45,45)
\qbezier(45,45)(65,65)(85,85)


\qbezier(5,5)(6.25,6.25)(7.5,7.5)
\qbezier(7.5,7.5)(8.75,10)(10,12.5)
\put(9,13.5){\makebox(0,0)[cc]{\footnotesize $D'$}}
\qbezier(10,12.5)(12.5,13.75)(15,15)
\qbezier(15,15)(50,50)(85,85)


\qbezier(5,5)(5.3125,5.3125)(5.625,5.625)
\qbezier(5.625,5.625)(5.9375,6.25)(6.25,6.875)
\qbezier(6.25,6.875)(6.875,7.1875)(7.5,7.5)
\qbezier(7.5,7.5)(46.25,46.25)(85,85)
\end{picture}

\caption{Graphs of the maps $a^{-k} ba^k$ (the graph
of $b$ is $OABCZ,$ the graph of $a^{-1}ba$ is $OCDEZ,$
the graph of $aba\inv$ is $OC'B'AZ,$ etc.)}
\end{figure}

Hence
$$
S_k=\supp(a^{-k}ba^{k})=\supp(x_0^{-2k} x_1 x_2^{-1} x_0^{2k}) =
\left(\frac 12 x_0^{2k}, \frac 12 x_0^{2k+2}\right)
$$
for all integer numbers $k.$ Note that
$$
\lim_{n \to +\infty} \frac 12 x_0^{2n}=1
$$
and that
$$
\lim_{n \to -\infty} \frac 12 x_0^{2n}=0.
$$
Then
\begin{equation} \label{tition}
[0,1] =\bigcup_{k \in \Z} \av S_k=\bigcup_{k \in \Z} \left[\frac 12 x_0^{2k}, \frac 12 x_0^{2k+2}\right],
\end{equation}
as required.

(ii) By (i), the maps $a^{-k} b a^k$ where $k$ runs over $\Z$ generate a free abelian group
of infinite rank. Hence every element $g$ of $G=\str{a,b}$ is uniquely written in
the form
\begin{equation} \label{els_o_G}
g=a^m \prod_{k \in \Z} a^{-k} b^{m_k} a^k
\end{equation}
where in the latter product only finitely many terms are non-trivial.
The element $b=x_1 x_0^{-1} x_1^{-1} x_0$ is in the commutator
subgroup $F'$ of $F.$ Then $g$ is congruent
to $a^m$ modulo $F'.$ As $F/F'=\str{x_0F',x_1F'}$ is a free abelian
group of rank two \cite{Belk,CFP}, $m \ne 0$ implies that
$g \ne 1.$ If $m=0,$ then $g=1$ if and only if all exponents $m_k$ are zeroes. Thus the generators
$a,b$ satisfy exactly the relations that
are corollaries of those given in the presentation of the group $\Z \wr \Z$ above.
\end{proof}

\begin{Prop} \label{G-is-def}
The group $G$ is first-order definable in $\PL_2([0,1])$
with the parameters $x_0,x_1.$
\end{Prop}

\begin{proof} We are going to prove that the subgroup $\str a$ generated
by $a$ and the subgroup $\str{a^{-k} b a^k : k \in \Z}$
are both definable in $\PL_2([0,1])$ by means of first-order
logic with the parameters $x_0$ and $x_1.$ This in view of
the formula \eqref{els_o_G} will imply that $G=\str{a,b}$
is $\{x_0,x_1\}$-definable, completing the proof.

Let $J=[\eta,\mu]$ be a closed interval with dyadic endpoints.
$\PL_2^{>}(J)$ (resp. $\PL_2^{\ge}(J)$) is the family $\{f\}$ of all maps in
$\PL_2(J)$ such that $\a f > \a$ (resp. $\a f \ge \a$) for all $\a \in (\eta,\mu).$
The following result is an immediate corollary of
Lemma 2.16 in \cite{Kass}.

\begin{Lem} \label{KassLem}
Let $z \in \PL_2^{>}(J).$ Then for every $q=2^m$ where
$m \in \Z$ there is at most one map in $\PL_2(J)$
commuting with $z$ whose initial slope is $q.$
\end{Lem}

By the initial slope of an $f$ in $\PL_2([\eta,\mu])$ we
mean the right-sided limit at $\eta$ of the derivative
of $f,$ or, less formally, the slope of the first
linear `component' of $f.$

Lemma \theLem\ implies that the centralizer of $x_0$ in $\PL_2([0,1])$
is the subgroup $\str{x_0}$ generated by $x_0$: indeed,
the initial slope of $x_0$ is $2$ and the initial
slope of a power $x_0^m$ of $x_0,$ obviously
commuting with $x_0,$ is $2^m.$ The powers
of $a=x_0^2$ are the squares of elements of $\str{x_0}.$ Thus
the family of powers of $a$ is first-order definable
in $\PL_2([0,1])$ with the parameter $x_0.$
It follows that the set $\{a^{-k} b a^k : k \in \Z\}$ is
definable in $F$ with the parameters $x_0$ and $x_1.$

We claim that
$$
\str{a^{-k} b a^k : k \in \Z}=C(\{a^{-k} b a^k : k \in \Z\})
$$
where $C(X)$ is the centralizer of a subset $X$ of $\PL_2([0,1]).$
If so, the subgroup $\str{a^{-k} b a^k : k \in \Z}$
is $\{x_0,x_1\}$-definable in $\PL_2([0,1]),$ being
the centralizer of an $\{x_0,x_1\}$-definable set.

\begin{Claim}
Let $k \in \Z.$ Then the map $g_k = a^{-k} b a^k|_{\avst S_k}$ where $S_k=\supp(a^{-k} b a^k)$
is in $\PL_2^{>}(\avst S_k)$ and its  intial slope
is $2.$
\end{Claim}

\begin{proof} The semi-group $\PL_2^{\ge}([0,1])$
is preserved under conjugation by elements of $\PL_2([0,1]).$
Consequently, $g|_S \in \PL_2^{>}(\av S),$ where $S=\supp(g),$
for any conjugate $g$ of $b,$ since $b=x_1x_2^{-1} \in \PL_2^{\ge}([0,1]).$

It is easily seen that both
maps $a=x_0^2$ and $a\inv=x_0^{-2}$ are linear
on the intervals $(0,1/8)$ and $(7/8,1)$ and hence differentiable
on these intervals. Suppose that $g \in \PL_2^{\ge}([0,1])$
is differentiable at $\a \in (0,1/8).$ Then, by the Chain Rule, the derivative
of $a g a\inv$ at the point $\a a\inv$ exists
and is equal to
\begin{equation}
\left.\frac{d a\inv}{dt}\right|_{t=\a g} \cdot
\left.\frac{d g}{dt}\right|_{t=\a} \cdot
\left.\frac{d a}{dt}\right|_{t=\a a\inv},
\end{equation}
that is, to the derivative of $g$ at $\a,$ provided that $\a g < 1/8$ and $\a a\inv < 1/8.$

By direct calculations we see that the second statement of
the Claim is true for maps $b,aba^{-1}, a^2ba^{-2}.$
The support $S_{-2}=(1/32,1/8)$ of the map $a^2ba^{-2}$ enters the interval
$(0,1/8).$ We prove the statement for all $k \le -2$ by induction: by (\theequation), for any
$k \le -2,$ the initial slope of the map $a^{-k} b a^k|_{S_k}$
is equal to the initial slope of the map $a^{-k+1} b a^{k-1}|_{S_{k-1}}$ .
The argument for positive $k$ is similar: we check
directly that the map $a^{-1} b a$ whose support
is contained in $(7/8,1)$ behaves as described,
and prove the statement for all $k \ge 1$ by induction.
\end{proof}

Now let $f \in C(\{a^{-k} b a^k : k \in \Z\}).$ Then for every
$k \in \Z$ the restriction $f_k$ of $f$ on $S_k$ commutes with $g_k.$
By Claim \theClaim\ and Lemma \ref{KassLem}, $f_k=g_k^{l_k}$ for a suitable
integer $l_k.$  Finally, part (i) of Lemma \ref{is_isom}  implies that
at most finitely many $l_k$ are non-zero, whence
$$
f=\prod_{k \in \Z} a^{-k} b^{l_k} a^k,
$$
and the result follows.
\end{proof}

\begin{Prop}
The elementary theory of Thompson's group $F$ is hereditarily
undecidable.
\end{Prop}

\begin{proof}
It is well-known that if a structure $\cM$ first-order interprets
with parameters a structure $\cN$ whose elementary theory
is hereditarily undecidable, then the elementary theory
of $\cM$ is too hereditarily undecidable \cite{Ersh}.

According to \cite{Nosk}, if a finitely generated almost solvable group $S$
is not almost abelian, then $S$ first-order interprets the
ring of integers $\Z$ (see also \cite{DeSi}.)
It is easy to see that the group $\Z \wr \Z$ is solvable;
in order to see that it is not almost abelian (that is,
does not have normal abelian subgroups of finite index)
one can work with the isomorphic group $G=\str{a,b}$ introduced above.
Suppose that $\varphi$ is a surjective homomorphism from $G$
onto a finite group of order $n.$ Then
$$
1=\varphi(a)^n=\varphi(b)^n = \varphi(a^n)=\varphi(b^n).
$$
The elements $a^n,b^n$ are not commuting, and hence the
kernel of $\varphi$ is not abelian.

Thus the elementary theory of the group $\Z \wr \Z$ is hereditarily
undecidable, for the elementary theory of
the ring of integers is hereditarily
undecidable \cite{Ersh,Ro}.

Now, by Proposition \ref{G-is-def}, Thompson's group $F$ first-order interprets
with parameters the group $\Z \wr \Z,$ and we are done.
\end{proof}

The authors would like to express their gratitude to Oleg
Belegradek for interesting discussions and to thank
Victor Guba and Mark Sapir for helpful information.


\begin{thebibliography}{9}


\bibitem{Belk} J.~Belk, Thompson's group $F$,
Ph. D. Thesis, Cornell Univ., 2004, available
at http://www.math.cornell.edu/~belk/Thesis.pdf.

\bibitem{CFP} J.~Cannon, W.~Floyd, and W.~Parry, Introductory notes to
Richard Thompson's groups, Enseign. Math. 42 (1996) 215--256.

\bibitem{Clea} S.~Cleary, Distortion of wreath products in some finitely presented groups,
2005, preprint, arXiv:math.GR/0505288.

\bibitem{DeSi} F.~Delon, P.~Simonetta, Undecidable wreath products and skew power series fields,
J. Symb. Logic 63 (1998) 237--246.

\bibitem{Ersh} Yu.~Ershov, Decision problems and constructivizable models, Nauka, Moscow, 1980.

\bibitem{GS} V.~Guba, M.~Sapir, On subgroups of the R.~Thompson group $F$ and other diagram groups,
Mat. Sb. 190 (1999) 3--60.

\bibitem{Kass} M.~Kassabov, F.~Matucci, Simultaneous conjugacy
problem in Thompson's group $F,$ 2006, preprint,
arXiv:math.GR/0607167.

\bibitem{Nosk} G.~Noskov, The elementary theory of a finitely generated almost solvable group,
Izv. Akad. Nauk SSSR Ser. Mat. 47 (1983) 498--517.

\bibitem{TGr40yrs} Problem list of the workshop `Thompson's Group at 40 Years' (Jan. 11-14 2004,
American Institute of Mathematics, Palo Alto, California), available at
http://aimath.org/WWN/thompsonsgroup/thompsonsgroup.pdf

\bibitem{Ro} R.~Robinson, Undecidable rings,
Trans. Amer. Math. Soc. 70 (1951) 137--159.
\end{thebibliography}
\end{document}